\documentclass[a4paper,12pt]{article}
\usepackage{amssymb}

\newcommand{\R}{\mathbb{R}}

\newcommand{\beq}{\begin{equation} }
\newcommand{\eqq}{\end{equation} }
\newcommand{\cuad}{{\sqcap\kern-.68em\sqcup}}

\newcommand{\equ}[1]{(\ref{#1})}

\newtheorem{definition}{Definition}[section]
\newtheorem{teo}{Theorem}[section]

\newtheorem{proposition}{Proposition}[section]

\newtheorem{lemma}{Lemma}[section]
\newtheorem{corollary}{Corollary}[section]
\newtheorem{remark}{Remark}[section]
\newcommand{\bremark}{\begin{remark} \em}
\newcommand{\eremark}{\end{remark} }

\def\beeq{\begin{equation}}
\def\eeq{\end{equation}}
\newcommand{\begeqaet}{\begin{eqnarray*}}
\newcommand{\eneqaet}{\end{eqnarray*}}

\begin{document}
\begin{center}{\bf  \large RADIAL SYMMETRY OF POSITIVE SOLUTIONS TO EQUATIONS INVOLVING THE FRACTIONAL LAPLACIAN }\medskip

\bigskip

\bigskip

{Patricio Felmer\ \  and\ \  Ying Wang
\\~\\}

Departamento de Ingenier\'{\i}a  Matem\'atica and Centro de
Modelamiento Matem\'atico,
UMR2071 CNRS-UChile,
 Universidad de Chile\\
{\sl   ( pfelmer@dim.uchile.cl   and   yingwang00@126.com ) }

\bigskip

\bigskip
\begin{abstract}
The aim of this paper is to study radial  symmetry and monotonicity properties  for positive solution of
elliptic equations involving the fractional Laplacian. We first consider the semi-linear Dirichlet problem
\begin{equation}\label{eq 021}
 (-\Delta)^{\alpha} u=f(u)+g,\ \ {\rm{in}}\ \ B_1, \quad u=0\ \ {\rm in}\ \ B_1^c,
\end{equation}
where $(-\Delta)^\alpha$ denotes the fractional Laplacian, $\alpha\in(0,1)$,  and
$B_1$ denotes the open unit ball centered at the origin in $\R^N$ with $N\ge2$.
The function $f:[0,\infty)\to\R$ is assumed to be locally Lipschitz continuous and  $g: B_1\to\R$ is  radially symmetric and decreasing in $|x|$.

In the second place we  consider radial symmetry of positive solutions  for  the equation
\begin{equation}\label{eq 021i}
 (-\Delta)^{\alpha} u=f(u),\ \ {\rm{in}}\ \ \R^N,
\end{equation}
with $u$ decaying  at infinity and  $f$
satisfying some extra hypothesis, but possibly being non-increasing.

Our third goal is to consider
radial  symmetry of positive solutions for system of the form
 \begin{equation}\label{eq 001}
\left\{ \arraycolsep=1pt
\begin{array}{lll}
(-\Delta)^{\alpha_1} u=f_1(v)+g_1,\ \ \ \ & {\rm{in}}\quad  B_1,\\[2mm]
(-\Delta)^{\alpha_2} v=f_2(u)+g_2,\ \ \ \ & {\rm{in}} \quad   B_1,\\[2mm]
 u=v =0,\ \ \ \   &  {\rm{in}}\quad   B_1^c,
\end{array}
\right.
\end{equation}
where $\alpha_1,\alpha_2\in(0,1)$, the functions $f_1$ and $f_2$ are locally Lipschitz continuous and increasing in $[0,\infty)$, and the functions $g_1$ and $g_2$ are radially symmetric and decreasing.

We prove our results
through the method of moving planes, using the recently proved ABP estimates for the fractional Laplacian. We
use a truncation technique to overcome the difficulty introduced by the non-local character of the differential operator in the application of the   moving planes.

\end{abstract}

\bigskip
\vspace{3mm}
\end{center}
\setcounter{equation}{0}
\section{Introduction}

The purpose of this paper is to study  symmetry and monotonicity properties of positive solutions  for
equations involving the fractional Laplacian  through the use of moving planes arguments. The first part of this article  is devoted to the following semi-linear Dirichlet problem
\begin{equation}\label{eq 1}
\left\{ \arraycolsep=1pt
\begin{array}{lll}
 (-\Delta)^{\alpha} u=f(u)+g,\ \ \ \ &
{\rm in} \quad B_1,\\[2mm]
 u=0,\quad & {\rm in}\quad B_1^c,
\end{array}
\right.
\end{equation}
where $B_1$ denotes the open unit ball centered at the origin in $\R^N,$ $N\ge2$ and
 $(-\Delta)^\alpha$ with  $\alpha\in(0,1)$ is the fractional Laplacian defined as
  \begin{equation}\label{DELTA}(-\Delta)^\alpha u(x)=P.V.\int_{\R^N}\frac{u(x)-u(y)}{|x-y|^{N+2\alpha}}dy,\end{equation}
$x\in B_1$. Here $P.V.$ denotes the principal value of the integral, that for notational simplicity we omit in what follows.

During the last  years, non-linear equations involving  general
integro-differential operators, especially, fractional Laplacian,
have been studied by many authors.  
Caffarelli and Silvestre
\cite{CS}
 gave a formulation of the fractional Laplacian through
Dirichlet-Neumann maps. Various regularity issues  for fractional elliptic
equations has been studied  by Cabr\'{e} and Sire \cite{CS0},  Caffarelli and Silvestre  \cite{CS1}, Capella, D\'{a}vila, Dupaigne and Sire  \cite{CDDS}, Ros-Oton and Serra  \cite{RS} and Silvestre  \cite{S}. Existence and related results were studied by  Cabr\'{e} and  Tan \cite{CT}, Dipierro, Palatucci and Valdinoci \cite{DPV}, Felmer, Quaas and Tan \cite{FQT}, and Servadei and Valdinoci \cite{SV}. 
Great attention has also been devoted to  symmetry results
for equations involving the fractional Laplacian in $\R^N$, such
as in the work by  Li \cite{L} and Chen, Li and Ou \cite{CLO1, CLO2},
where the method of moving planes in integral form has been
developed to treat various equations and systems, see
also Ma and Chen \cite{MC}. On the other hand,  using the local formulation of Caffarelli and Silvestre, Cabr\'{e} and  Sire \cite{CS5} applied
the sliding method to obtain  symmetry results for nonlinear
equations with fractional laplacian and 
 Sire and Valdinoci \cite{S1} 
 studied symmetry properties  for a boundary reaction problem via a
geometric inequality.
 Finally,  in
\cite{FQT}  the authors used the method of moving planes in integral form to
prove  symmetry results for
\begin{equation}\label{eq 1wy}
(-\Delta)^\alpha u+u=h(u)\ \ \rm{in}\ \ \R^N,
\end{equation}
taking advantage of  the representation formula for $u$ given by
 $$u(x)=\mathcal{K}\ast h(u)(x),\quad x\in\R^N,$$
 where   the kernel $\mathcal{K}$,  associated to the linear part of the equation,   plays a key role
 in the  arguments.
This approach is not possible to be used for problem (\ref{eq 1}),
since a similar representation formula is not available in general.

The study of radial symmetry and monotonicity of positive solutions for non-linear elliptic equations in bounded domains
 using the moving planes method based on the Maximum Principle was initiated with the work by Serrin \cite{JS} and Gidas, Ni and Nirenberg \cite{GNN},
 with important subsequent advances by  Berestycki and Nirenberg \cite{BN}. We refer to the review by
Pacella and Ramaswamy \cite{PR}  for a more complete discussion of
the method and it various applications. In \cite{BN} the Maximum
Principle for small domain, based on the  Aleksandrov-Bakelman-Pucci
(ABP) estimate, was used as a tool to obtain much general results,
specially avoiding regularity hypothesis on the domain.
 In a recent article Guillen and Schwab,  \cite{GS}, provided an ABP estimate   for a  class of fully non-linear elliptic integro-differential equations.
Motivated by this work, we obtain a version of the Maximum Principle for small domain and we apply it with
the  moving planes method  as in \cite{BN} to prove symmetry and monotonicity properties for
 positive solutions to problem (\ref{eq 1}) in the ball and in more general domains.

We consider the following hypotheses on the functions $f$ and $g$:
\begin{itemize}
\item[$(F1)\ $]
The function  $f:[0,\infty)\to\R$ is locally Lipschitz.
\item[$(G)\ $]
The function  $g: {B_1}\to\R$ is  radially symmetric and decreasing in $|x|$.
\end{itemize}

Before stating our first theorem we make precise the notion of solution that we use in this article. We say that a continuous function $u:\R^N\to\R$ is a classical solution of equation (\ref{eq 1}) if the fractional Laplacian of $u$ is defined at any point of $B_1$, according to the definition given in \equ{DELTA}, and if $u$ satisfies  the equation and the external condition in a pointwise sense. This notion of solution is extended in a natural way to the other equations considered in this paper.

Now we are ready for our first theorem on radial symmetry and
monotonicity properties for positive solutions of equation
 (\ref{eq 1}). It states   as
follows:
\begin{teo}\label{teo 1}
 Assume that
the functions $f$ and $g$ satisfy $(F1)$ and $(G)$, respectively. If  $u$ is a positive classical solution of (\ref{eq 1}), then
$u$ must be radially symmetric and strictly decreasing in $r=|x|$
for $r\in(0,1)$.
\end{teo}

The proof of Theorem \ref{teo 1} is given in  Section \S 3, where we prove a more general symmetry and monotonicity result for equation (\ref{eq 1}) on a general domain $\Omega$, which is convex and symmetric in one direction, see Theorem \ref{remark 111}.

We devote the  second part of this article to study  symmetry
results for a non-linear equation as (\ref{eq 1}), but in $\R^N$ and
with $g\equiv0$. For the problem in $\R^N$,  the moving planes
procedure has to start  a different way because we cannot use the
Maximum Principle for small domain. We  refer to the work by
 Gidas, Ni and Nirenberg \cite{GNN2}, Berestycki and Nirenberg \cite{BN},
Li \cite{LCM}, and Li and Ni \cite{LN}, where these results were studied assuming some additional hypothesis on $f$,
 allowing for decay properties of the solution $u$.
A general result in this direction was obtained by Li \cite{LCM} for the equation
 $$-\Delta u=f(u),\ \ \ \mbox{in} \ \ \R^N,$$
with $u$ decaying  at infinity and $f$ satisfying the following hypothesis:
\begin{itemize}
\item[$(F2)\ $]
\begin{enumerate}\item[]
 There exists $s_0>0,\ \gamma>0$ and $C>0$ such that
  \begin{equation}\label{eq rq2}
 \frac{f(v)-f(u)}{v-u}\leq C(u+v)^\gamma\ \ \ \ \mbox{for all } \quad 0<u<v<s_0.
 \end{equation}
\end{enumerate}
\end{itemize}
Motivated by these results, we are interested in similar  properties
of positive solutions for equations involving the fractional
Laplacian under assumption (F2). Here is our second main theorem.

\begin{teo}\label{teo r1}
Assume that $\alpha\in(0,1),$ $N\geq2$, the function $f$ satisfies $(F1)-(F2)$ and $u$  is
a positive classical solution for  the equation
\begin{equation}\label{eq wyq1}
\left\{ \arraycolsep=1pt
\begin{array}{lll}
 (-\Delta)^{\alpha} u=f(u)\ \ \ in\ \ \R^N,\\[2mm]
 u>0\ \ in\ \ \R^N,\ \ \lim_{|x|\to\infty}u(x)=0.
\end{array}
\right.
\end{equation}
Assume further that there exists
\begin{equation}\label{m}m>\max\{\frac{2\alpha}{\gamma},\frac{N}{\gamma+2}\}
\end{equation} such that  $u$ satisfies
\begin{equation}\label{ying1.4}
u(x)=O(\frac{1}{|x|^m}), \quad\mbox{as}\quad |x|\to\infty,
\end{equation}
then $u$ is radially symmetric and strictly decreasing about some
point in $\R^N$.
\end{teo}

 In \cite{FQT}, Felmer, Quaas and Tan studied symmetry of positive solutions using
the integral form of the moving planes method, assuming that the
function $f$ is such that $h(\xi)\equiv f(\xi)+\xi$  is
super-linear, with sub-critical growth at infinity and
\begin{itemize}
\item[$(H)\ $]
\begin{enumerate}\item[]
$h\in C^1(\R),$ increasing and there exists $\tau>0$ such that
 $$\lim_{v\to0}\frac{h'(v)}{v^{\tau}}=0.$$
\end{enumerate}
\end{itemize}
We see that Theorem \ref{teo r1} generalizes Theorem 1.3  in
\cite{FQT},  since here we do not assume $f$ is differentiable and
we do not require $h$  to be increasing. In Section \S 4 we  present
an extension of Theorem \ref{teo r1} to  $f(\xi)=\xi^p-\xi^q$, with
$0<q<1<p$, that is not    covered by the results in  \cite{FQT}
either,  see Theorem  \ref{teo 42}. This non-linearity was studied
by Valdebenito in \cite{valdebenito}, where decay and symmetry
results were obtained using local extension as in Caffarelli and
Silvestre \cite{CS} and regular moving planes.

For the  particular case $f(u)=u^p$, for some $p>1$, we see that (H) is not satisfied, but that (F2) does hold. Thus, if we knew  the solution of \equ{eq wyq1} satisfies  decay assumption   \equ{ying1.4}  in this setting, we would have symmetry results in these cases. See \cite{GNN2} and \cite{LCM} for the proof of decay properties in the case of the Laplacian.


The third part  of this paper is devoted to investigate the radial symmetry of  non-negative solutions for the following
system of non-linear equations involving  fractional Laplacians
 with different orders,
\begin{equation}\label{eq 21}
\left\{ \arraycolsep=1pt
\begin{array}{lll}
(-\Delta)^{\alpha_1} u=f_1(v)+g_1,\ \ \ \ &\mbox{in}\quad
 B_1,\\[2mm]
(-\Delta)^{\alpha_2} v=f_2(u)+g_2,\ \ \ \ &\mbox{in}\quad
 B_1,\\[2mm]
 u=v=0,&\mbox{in}\quad  B_1^c,
\end{array}
\right.
\end{equation}
where  $\alpha_1,\alpha_2\in(0,1)$.
We have following  results for system (\ref{eq 21}):
\begin{teo}\label{teo 2}
 Suppose $f_1$ and $f_2$ are locally Lipschitz continuous and increasing functions defined in $[0,\infty)$ and $g_1$ and $g_2$ satisfy (G). Assume that  $(u,v)$ are positive, classical solutions  of system (\ref{eq 21}), then $u$ and $v$ are
radially symmetric and strictly decreasing in  $r=|x|$ for
$r\in(0,1)$.
\end{teo}

We prove  our theorems using the moving planes method. The main
difficulty comes from the fact that the fractional Laplacian is a
non-local operator, and consequently  Maximum Principle and
Comparison Results require information on the solutions in the whole
complement of the domain, not only at the boundary.  To overcome
this difficulty, we introduce a new  truncation technique which is
well adapted to be used with the   method of moving planes.

The rest of the paper is organized as follows. In Section \S 2, we
recall  the ABP estimate for equations involving fractional Laplacian, as proved in \cite{GS} and
we prove a form of Maximum Principle for domains with small measure.
In Section \S3, we prove Theorem \ref{teo 1} by the moving planes method and we extend our symmetry results to general domains with one dimensional convexity and symmetry properties.
In Section \S4, the radial symmetry of solutions for equation (\ref{eq wyq1})  in $\R^N$ is obtained. An extension to a non-lipschitzian non-linearity is given.
In Section \S5, we complete the proof of  Theorem \ref{teo 2}.
And finally, Section \S 6  is devoted to discuss (\ref{eq 1})
for a non-local operator with non-homogeneous kernel.

 \setcounter{equation}{0}
\section{Preliminaries}
A key tool in the use of  the  moving planes method is the Maximum Principle for small domain, which is a consequence of the ABP estimate. In
\cite{GS}, Guillen and  Schwab showed an ABP estimate (see
Theorem 9.1) for general integro-differential operators.   In this section we recall this estimate in the case of the
fractional Laplacian in any open and bounded domain. Then we obtain the  Maximum Principle for small domains.

We start with the  ABP estimate for the fractional Laplacian, which is stated as follows:
\begin{proposition}\label{pro abp1}
Let $\Omega$ be a bounded, connected open subset of $\R^N$. Suppose that $h:\Omega\to\R$ is in $ L^\infty(\Omega)$ and $w\in L^\infty(\R^N) $ is a classical solution of
\begin{equation}\label{abp 1}
\left\{ \arraycolsep=1pt
\begin{array}{lll}
\Delta^\alpha w(x)\le h(x),\ \ \ \ &
x\in \Omega,\\[2mm]
 w(x)\ge 0,& x\in \R^N\setminus \Omega.
\end{array}
\right.
\end{equation}

 Then there exists a positive constant $C$, depending on  $N$ and $\alpha$, such that
 \begin{equation}\label{abp}
 -\inf_{\Omega}w\le Cd^\alpha\|h^+\|_{L^\infty(\Omega)}^{1-\alpha}\|h^+\|_{L^{N}(\Omega)}^\alpha,
 \end{equation}
 where $d=\mbox{diam}({\Omega})$ is the diameter of $\Omega$ and  $h^+(x)=\max\{h(x),0\}$.

 Here and in what follows we write $\Delta^\alpha w(x)=-(-\Delta)^\alpha w(x).$
 \end{proposition}

We have the following corollary
\begin{corollary}\label{cor 1}
Under the assumptions of  Proposition \ref{pro abp1}, with $\Omega$ not necessarily connected, we have
 \begin{equation}\label{abpx}
 -\inf_{\Omega}w\le C d^{\alpha}\|h^+\|_{L^\infty(\Omega)} |\Omega|^{\frac\alpha N}.
 \end{equation}
\end{corollary}
\noindent{\bf Proof.}
We let
 $w_0\in L^\infty(\R^N)$ be a classical solution of
 \begin{equation}\label{abp 2}
\left\{ \arraycolsep=1pt
\begin{array}{lll}
\Delta^\alpha w_0(x)= \|h^+\|_{L^\infty(\Omega)}\chi_{\Omega}(x),\ \ \ \ &
x\in B_d(x_0),\\[2mm]
 w_0(x)= 0,& x\in \R^N\setminus B_d(x_0),
\end{array}
\right.
\end{equation}
where $x_0\in\Omega$ and $\Omega\subset B_d(x_0)$. We observe that $B_d(x_0)$ is connected and that $w_0\le 0$ in all $\R^N$.
By Comparison Principle, we see that
$$\inf_{\R^N} w_0\le \inf_{\R^N} w,$$
where $w$ is the solution of (\ref{abp 1}).
Then we use  Proposition \ref{pro abp1} to obtain that
$$-\inf_{\R^N} w_0=-\inf_{B_d(x_0)}w_0\le C (2d)^{\alpha}\|h^+\|_{L^\infty(\Omega)} |\Omega|^{\frac\alpha N}$$
and then we conclude
$$-\inf_{\Omega}w=-\inf_{\R^N}w\le C d^{\alpha}\|h^+\|_{L^\infty(\Omega)} |\Omega|^{\frac\alpha N}.\qquad\Box$$

\begin{remark}\label{remark abp1}
We notice that, under a possibly different constant $C>0$,  the ABP estimate for problem (\ref{abp 1}) with $\alpha=1$
$$\left\{ \arraycolsep=1pt
\begin{array}{lll}
\Delta w(x)\le h(x),\ \ \ \ &
x\in \Omega,\\[2mm]
 w(x)\ge 0,& x\in \partial \Omega,
\end{array}
\right.$$
 is
precisely (\ref{abp}) with $\alpha=1$.
 \end{remark}

 As a consequence of the ABP estimate just recalled, we have the
Maximum Principle for small domain, which is stated as follows:
\begin{proposition}\label{teo 3} Let $\Omega$  be an
open and bounded subset of $\R^N$. Suppose that
$\varphi:\Omega\to\R$ is in $L^\infty(\Omega)$ and
$w\in L^\infty(\R^N)$ is a classical solution of
\begin{equation}\label{eq a1}
\left\{ \arraycolsep=1pt
\begin{array}{lll}
\Delta^\alpha  w(x)\le \varphi(x)w(x),\ \ \ \ &
x\in \Omega,\\[2mm]
 w(x)\ge 0,& x\in \R^N\setminus \Omega.
\end{array}
\right.
\end{equation}

Then there is $\delta>0$ such that whenever $|\Omega^-|\le \delta$,
$w$ has to be  non-negative in $\Omega$. Here $\Omega^-=\{x\in\Omega\ | \ w(x)<0\}$.
 \end{proposition}
 \noindent{\bf Proof.}
By (\ref{eq a1}), we observe that
 $$\left\{ \arraycolsep=1pt
 \begin{array}{lll}
\Delta^\alpha  w(x)\le \varphi(x)w(x),\ \ \ \ &
x\in \Omega^-,\\[2mm]
 w(x)\ge 0,& x\in \R^N\setminus \Omega^-.
\end{array}
\right.
$$
Then, using Corollary \ref{cor 1} with $h(x)=\varphi(x)w(x)$, we
obtain that
 \begin{eqnarray*}
\|w\|_{L^\infty(\Omega^-)} =-\inf_{\Omega^-}w &\le& Cd_0^{\alpha}
\|(\varphi w)^+\|_{L^\infty(\Omega^-)} |\Omega^-|^{\frac\alpha N},
\end{eqnarray*}
where constant $C>0$ depends on  $N$ and $\alpha$. Here  $d_0=\mbox{diam}(
\Omega^-)$. Thus
$$
\|w\|_{L^\infty(\Omega^-)}\le
Cd_0^{\alpha}\|\varphi\|_{L^\infty(\Omega)}\|w\|_{L^\infty(\Omega^-)}|\Omega^-|^{\frac\alpha
N}.
$$
We see that,  if $|\Omega^-|$ is  such that
$Cd_0^{\alpha}\|\varphi\|_{L^\infty(\Omega)}|\Omega^-|^{\alpha/N}<1$,
then we must have that
$$\|w\|_{L^\infty(\Omega^-)}=0.$$
This implies that  $|\Omega^-|=0$ and
since $\Omega^-$ is open, we have   $\Omega^-=\emptyset$, completing the
 proof.  \hfill$\Box$\\

\setcounter{equation}{0}
\section{Proof of Theorem \ref{teo 1}.}

In this section we provide  a proof of Theorem \ref{teo 1} on the radial symmetry and monotonicity of positive solutions to equation (\ref{eq 1}) in the unit ball.  For this purpose we use the   of moving planes method, for which we give some preliminary notation. We define
\begin{equation}\label{d1}\Sigma_\lambda=\{x=(x_1,x')\in B_1\  |\  x_1>\lambda\},\end{equation}
\begin{equation}\label{d2}T_\lambda=\{x=(x_1,x')\in \R^{N}\ |\  x_1=\lambda\},\end{equation}
\begin{equation}\label{d3}u_\lambda(x)=u(x_\lambda) \quad\mbox{and}\quad  w_\lambda(x)=u_\lambda(x)-u(x),\end{equation}
where $\lambda\in (0,1)$ and $x_\lambda=(2\lambda-x_1,x')$  for
$x=(x_1,x')\in\R^N.$ For any subset $A$ of $\R^N$, we write
$A_\lambda=\{x_\lambda:\, x\in A\}$, the reflection of $A$ with
regard to $T_\lambda$.

\medskip

\noindent{\bf Proof of  Theorem \ref{teo 1}.}
 We divide the proof in three steps.

\noindent\textbf{Step 1:} We prove that if $\lambda\in (0,1)$ is close to $1$, then
$w_\lambda>0$ in $\Sigma_\lambda$. For this purpose, we start proving that
 if $\lambda\in (0,1)$ is close to 1, then $w_\lambda\ge0$ in $\Sigma_\lambda$.
If we define $\Sigma_\lambda^-=\{x\in \Sigma_\lambda\ |\ w_\lambda(x)<0\},$
then we just need to prove that if
$\lambda\in (0,1)$ is close to 1 then
\begin{equation}\label{eq 11}
\Sigma_\lambda^-=\emptyset.
\end{equation}

By contradiction, we assume (\ref{eq 11}) is not true, that is
$\Sigma^-_\lambda\not=\emptyset$.
 We denote
\begin{equation}\label{eq c}
w_\lambda^+(x)=\left\{ \arraycolsep=1pt
\begin{array}{lll}
w_\lambda(x),\ \ \ \ &
x\in \Sigma_\lambda^-,\\[2mm]
0,& x\in \R^N\setminus \Sigma_\lambda^-,
\end{array}
\right.
\end{equation}
\begin{equation}\label{eq 4.01}
w_\lambda^-(x)=\left\{ \arraycolsep=1pt
\begin{array}{lll}
0,\ \ \ \ &
x\in \Sigma_\lambda^-,\\[2mm]
w_\lambda(x),\ \ \ \ & x\in \R^N\setminus \Sigma_\lambda^-
\end{array}
\right.
\end{equation}
and we observe that
 $w_\lambda^+(x)=w_\lambda(x)-w_\lambda^-(x)$ for all $x\in\R^N.$
Next we claim that for all $0<\lambda<1$, we have
\begin{equation}\label{claim1}
 (-\Delta)^{\alpha} w_\lambda^-(x)\le0,\ \ \ \ \forall\ x\in \Sigma_\lambda^-.
 \end{equation}
By direct computation, for
 $x\in \Sigma_\lambda^-$, we have
 \begin{eqnarray*}
 (-\Delta)^{\alpha} w_\lambda^-(x)
&=&\int_{\R^N}\frac{w_\lambda^-(x)-w_\lambda^-(z)}{|x-z|^{N+2\alpha}}dz
=-\int_{\R^N\setminus\Sigma_\lambda^-}\frac{w_\lambda(z)}{|x-z|^{N+2\alpha}}dz
\\&=&-\int_{(B_1\setminus(B_1)_\lambda) \cup ((B_1)_\lambda\setminus B_1)}\frac{w_\lambda(z)}{|x-z|^{N+2\alpha}}dz
\\&&-\int_{(\Sigma_\lambda\setminus\Sigma_\lambda^-) \cup (\Sigma_\lambda\setminus\Sigma_\lambda^-)_\lambda}\frac{w_\lambda(z)}{|x-z|^{N+2\alpha}}dz
-\int_{(\Sigma_\lambda^-)_\lambda}\frac{w_\lambda(z)}{|x-z|^{N+2\alpha}}dz
\\&=&-I_1-I_2-I_3.
\end{eqnarray*}
We look at each of these integrals separately. Since $u=0\ in \ (B_1)_\lambda\setminus B_1$ and
$u_\lambda=0\ in \ B_1\setminus (B_1)_\lambda$,
we have
\begin{eqnarray*}
I_1
&=&\int_{(B_1\setminus(B_1)_\lambda) \cup ((B_1)_\lambda\setminus B_1)}\frac{w_\lambda(z)}{|x-z|^{N+2\alpha}}dz
\\&=&\int_{(B_1)_\lambda\setminus B_1}\frac{u_\lambda(z)}{|x-z|^{N+2\alpha}}dz-\int_{ B_1\setminus(B_1)_\lambda}\frac{u(z)}{|x-z|^{N+2\alpha}}dz
\\&=&\int_{(B_1)_\lambda\setminus B_1}u_\lambda(z)(\frac{1}{|x-z|^{N+2\alpha}}-\frac{1}{|x-z_\lambda|^{N+2\alpha}}))dz\ge 0,
\end{eqnarray*}
since $u_\lambda\geq0$ and $|x-z_\lambda|>|x-z|$ for all  $x\in \Sigma_\lambda^-$ and $z\in (B_1)_\lambda\setminus B_1.$
In order to study the sign of $I_2$ we first observe that $w_\lambda(z_\lambda)=-w_\lambda(z)$ for any $z\in\R^N$. Then
\begin{eqnarray*}
I_2
&=&\int_{(\Sigma_\lambda\setminus\Sigma_\lambda^-) \cup(\Sigma_\lambda\setminus\Sigma_\lambda^-)_\lambda}\frac{w_\lambda(z)}{|x-z|^{N+2\alpha}}dz
\\&=&\int_{\Sigma_\lambda\setminus\Sigma_\lambda^-}\frac{w_\lambda(z)}{|x-z|^{N+2\alpha}}dz
+\int_{\Sigma_\lambda\setminus\Sigma_\lambda^-}\frac{w_\lambda(z_\lambda)}{|x-z_\lambda|^{N+2\alpha}}dz
\\&=&\int_{\Sigma_\lambda\setminus\Sigma_\lambda^-}w_\lambda(z)(\frac{1}{|x-z|^{N+2\alpha}}-\frac{1}{|x-z_\lambda|^{N+2\alpha}})dz\ge 0,
\end{eqnarray*}
since  $w_\lambda\ge0$ in $\Sigma_\lambda\setminus\Sigma_\lambda^-$ and $|x-z_\lambda|>|x-z|$ for all $x\in \Sigma_\lambda^-$ and $z\in
\Sigma_\lambda\setminus\Sigma_\lambda^-.$
Finally, since
  $w_\lambda(z)<0$ for $z\in \Sigma_\lambda^-$, we have
\begin{eqnarray*}
I_3
&=&\int_{(\Sigma_\lambda^-)_\lambda}\frac{w_\lambda(z)}{|x-z|^{N+2\alpha}}dz
=\int_{\Sigma_\lambda^-}\frac{w_\lambda(z_\lambda)}{|x-z_\lambda|^{N+2\alpha}}dz
\\&=&-\int_{\Sigma_\lambda^-}\frac{w_\lambda(z)}{|x-z_\lambda|^{N+2\alpha}}dz
\ge0.
\end{eqnarray*}
Hence, we obtain \equ{claim1}, proving the claim.
Now we apply  \equ{claim1} and linearity of the fractional Laplacian to   obtain that, for $ x\in\Sigma_\lambda^-,$
\begin{equation}\label{eq e}
(-\Delta)^{\alpha} w_\lambda^+(x)\ge (-\Delta)^{\alpha} w_\lambda(x)
=(-\Delta)^{\alpha} u_\lambda(x)-(-\Delta)^{\alpha}  u(x).
\end{equation}
Combining  equation (\ref{eq 1}) with (\ref{eq e}) and  (\ref{eq c}), for $x\in\Sigma_\lambda^-$ we have
\begin{eqnarray*}
(-\Delta)^{\alpha} w_\lambda^+(x) &\geq&(-\Delta)^{\alpha}
u_\lambda(x)-(-\Delta)^{\alpha} u(x)
\\&=& f(u_\lambda(x))+g(x_\lambda)-f(u(x))-g(x)
\\&=&\frac{f(u_\lambda(x))-f(u(x))}{u_\lambda(x)-u(x)}w_\lambda^+(x)+g(x_\lambda)-g(x).
\end{eqnarray*}
Let us define
 $\varphi(x)=-({f(u_\lambda(x))-f(u(x))})/({u_\lambda(x)-u(x)})$  for $x\in\Sigma_\lambda^-$. By assumption $(F1)$, we have that $\varphi\in
L^\infty(\Sigma_\lambda^-)$. By assumption $(G)$, we have that  $g(x_\lambda)\geq g(x)$, since
for all $x\in\Sigma_\lambda^-$ and
$0<\lambda<1$, we have $|x|>|x_\lambda|$.
Hence, we have
\begin{equation}\label{eq f}
\Delta^{\alpha}  w_\lambda^+(x)\le \varphi(x)w_\lambda^+(x),\ \
x\in\Sigma_\lambda^-
\end{equation}
and  since
$w_\lambda^+=0$ in $(\Sigma_\lambda^-)^c$
we may apply Proposition \ref{teo 3}. Choosing $\lambda\in (0,1)$ close enough to $1$ we find that
$|\Sigma_\lambda^-|$ is small and then
$$w_\lambda=w_\lambda^+\geq0\ \ \ \ \mbox{in} \ \ \Sigma_\lambda^-.$$
But this is a contradiction with our assumption so we have
$$w_\lambda\geq0\ \ \ in\ \ \Sigma_\lambda.$$

In order to complete Step 1, we claim that for $0<\lambda<1$, if
 $w_\lambda\ge0$ and $w_\lambda\not\equiv0$
 in $\Sigma_\lambda$, then  $w_\lambda>0$
 in $\Sigma_\lambda$. Assuming the claim is true, we complete the proof, since the function $u$ is positive in $B_1$ and
 $u=0$ on $\partial B_1$, so that $w_\lambda$ is positive in $\partial B_1\cap \partial \Sigma_\lambda$ and then, by continuity   $w_\lambda\not =0$ in $\Sigma_\lambda$.


Now we prove the claim.   Assume there exists $x_0\in \Sigma_\lambda$ such that
$w_\lambda(x_0)=0,$ that is,  $u_\lambda(x_0)=u(x_0)$. Then we have that
\begin{eqnarray*}
(-\Delta)^{\alpha} w_\lambda(x_0)&=& (-\Delta)^{\alpha}
u_\lambda(x_0)-(-\Delta)^{\alpha} u(x_0)
=g((x_0)_\lambda)-g(x_0).
\end{eqnarray*}
Since $x_0\in \Sigma_\lambda$, we have   $|x_0|>|(x_0)_\lambda|$, then
by assumption $(G)$ we have $g((x_0)_\lambda)\geq g(x_0)$ and thus
\begin{equation}\label{eq 7}
(-\Delta)^{\alpha} w_\lambda(x_0)\geq0.
\end{equation}
On the other hand, defining
$A_\lambda=\{(x_1,x')\in\R^N\ |\  x_1>\lambda\}$,
since $w_\lambda(z_\lambda)=-w_\lambda(z)$ for any $z\in\R^N$ and $w_\lambda(x_0)=0$, we find
\begin{eqnarray*}
(-\Delta)^{\alpha} w_\lambda(x_0)
&=&-\int_{A_\lambda}\frac{w_\lambda(z)}{|x_0-z|^{N+2\alpha}}dz-\int_{\R^N\setminus A_\lambda}\frac{w_\lambda(z)}{|x_0-z|^{N+2\alpha}}dz
\\&=&-\int_{A_\lambda}\frac{w_\lambda(z)}{|x_0-z|^{N+2\alpha}}dz-\int_{A_\lambda}\frac{w_\lambda(z_\lambda)}{|x_0-z_\lambda|^{N+2\alpha}}dz
\\&=&-\int_{A_\lambda}w_\lambda(z)(\frac{1}{|x_0-z|^{N+2\alpha}}-\frac{1}{|x_0-z_\lambda|^{N+2\alpha}})dz.
\end{eqnarray*}
Since $|x_0-z_\lambda|>|x_0-z|$ for $z\in A_\lambda$ ,
$w_\lambda(z)\ge0$ and $w_\lambda(z)\not\equiv0$ in $A_\lambda$, from here we get
\begin{equation}\label{ff}
(-\Delta)^{\alpha} w_\lambda(x_0)<0,\end{equation}
which contradicts  (\ref{eq 7}), completing the proof of the claim.

\noindent\textbf{Step 2:} We define  $\lambda_0=\inf\{\lambda\in(0,1)\ |\  w_\lambda>0\ \ \rm{in}\ \ \Sigma_\lambda\}$ and we prove that $\lambda_0=0$. Proceeding by contradiction, we assume that $\lambda_0>0$, then
 $w_{\lambda_0}\ge0$ in $\Sigma_{\lambda_0}$ and
$w_{\lambda_0}\not\equiv0$ in $\Sigma_{\lambda_0}$. Thus, by the claim just proved above,
 we have $w_{\lambda_0}>0$ in $\Sigma_{\lambda_0}$.

Next we claim that if $w_\lambda>0$ in $\Sigma_\lambda$ for
$\lambda\in(0,1)$, then there exists $\epsilon\in(0,\lambda)$ such
that $w_{\lambda_\epsilon}>0$ in $\Sigma_{\lambda_\epsilon}$, where
$\lambda_\epsilon=\lambda-\epsilon$. This claim directly implies that $\lambda_0=0$, completing Step 2.


Now we  prove the claim.
Let $D_\mu=\{x\in\Sigma_\lambda\ | \ dist(x,\partial\Sigma_\lambda)\ge \mu\}$ for $\mu>0$ small. Since
$w_\lambda>0$ in $\Sigma_\lambda$ and $D_\mu$ is compact, then there
exists $\mu_0>0$ such that $w_\lambda\ge \mu_0$ in $D_\mu$. By
continuity of $w_\lambda(x)$, for $\epsilon>0$ small
enough and  denoting $\lambda_\epsilon=\lambda-\epsilon,$ we have that
$$w_{\lambda_\epsilon}(x)\ge0\ \ \rm{in}\ \ D_\mu.$$ As a
consequence, $$\Sigma_{\lambda_\epsilon}^-\subset
\Sigma_{\lambda_\epsilon}\setminus D_\mu$$ and
$|\Sigma_{\lambda_\epsilon}^-|$ is small if $\epsilon$ and $\mu$ are small.
Using \equ{claim1} and proceeding as in  Step 1, we have for all $x\in \Sigma_{\lambda_\epsilon}^-$
that \begin{eqnarray*}
(-\Delta)^{\alpha} w_{\lambda_\epsilon}^+(x) &=&(-\Delta)^{\alpha}
u_{\lambda_\epsilon}(x)-(-\Delta)^{\alpha}  u(x)-(-\Delta)^{\alpha}
w_{\lambda_\epsilon}^-(x)
\\&\ge& (-\Delta)^{\alpha} u_{\lambda_\epsilon}(x)-(-\Delta)^{\alpha} u(x)
\\&=&\varphi(x)w_{\lambda_\epsilon}^+(x)+g(x_\lambda)-g(x)
\ge \varphi(x)w_{\lambda_\epsilon}^+(x),
\end{eqnarray*}
where
$\varphi(x)=\frac{f(u_{\lambda_\epsilon}(x))-f(u(x))}{u_{\lambda_\epsilon}(x)-u(x)}$ is bounded  by assumption $(F1)$.

Since $w_{\lambda_\epsilon}^+=0$ in
$(\Sigma_{\lambda_\epsilon}^-)^c$ and
$|\Sigma_{\lambda_\epsilon}^-|$ is small, for  $\epsilon$ and $\mu$
small,  Proposition \ref{teo 3}  implies that  $w_{\lambda_\epsilon}\ge0$ in
$\Sigma_{\lambda_\epsilon}$.  Thus, since $\lambda_\epsilon>0$ and
$w_{\lambda_\epsilon}\not\equiv0$ in $\Sigma_{\lambda_\epsilon}$, as before we have
 $w_{\lambda_\epsilon}>0$
 in $\Sigma_{\lambda_\epsilon}$, completing the proof of the claim.

\noindent\textbf{Step 3:} By Step 2, we have $\lambda_0=0$, which
implies that $u(-x_1,x')\ge u(x_1,x')$ for $x_1\ge0.$
Using the same argument from the other side, we conclude that $u(-x_1,x')\le u(x_1,x')$ for $ x_1\ge0$ and then
$u(-x_1,x')= u(x_1,x')$ for $x_1\ge0.$ Repeating this procedure in all directions we obtain radial symmetry of $u$.

Finally, we prove $u(r)$ is strictly decreasing in  $r\in (0,1)$. Let us consider  $0<x_1<\widetilde{x}_1<1$
and let  $\lambda=\frac{x_1+\widetilde{x}_1}{2}$. Then, as proved above   we have
$$w_\lambda(x)>0\ \ \mbox{for}\ \ x\in\Sigma_{\lambda}.$$
Then
\begin{eqnarray*}
0<w_\lambda(\widetilde{x}_1,0,\cdots,0)
&=&u_\lambda(\widetilde{x}_1,0,\cdots,0)-u(\widetilde{x}_1,0,\cdots,0)
\\&=&u(x_1,0,\cdots,0)-u(\widetilde{x}_1,0,\cdots,0),
\end{eqnarray*}
that is
$u(x_1,0,\cdots,0)>u(\widetilde{x}_1,0,\cdots,0).$ Using the radial symmetry of $u$, we conclude from here the monotonicty of $u$.
 \hfill$\Box$\\

The proof of  Theorem \ref{teo 1} can be applied directly to prove  symmetry results
 for problem (\ref{eq 1}) in more general domains.
We have the following  definition
 \begin{definition}\label{def 111}
 We say that domain $\Omega\subset\R^N $ is convex in the $x_1$ direction:\\
$$ (x_1,x'),(x_1,y')\in\Omega\Rightarrow
(x_1,tx'+(1-t)y')\in\Omega,\ \ \forall\ t\in(0,1).$$
 \end{definition}
Now we state the more general theorem:
 \begin{teo}\label{remark 111}
Let  $\Omega\subset\R^N (N\ge2)$ is an open and bounded set. Assume further that $\Omega$ is convex in the $x_1$ direction and symmetric with respect to
the plane $x_1=0$.
 Assume that
the function $f$ satisfies $(F1)$ and  $g$ satisfies
\begin{itemize}
\item[$(\widetilde{G})\ $]
 The function
 $g: \Omega\to\R$ is  symmetric with respect to $x_1=0$ and decreasing in the $x_1$ direction,
for $x=(x_1,x')\in\Omega$, $x_1>0$.
\end{itemize}
Let $u$ be a positive classical solution  of
\begin{equation}\label{eq omega}
\left\{ \arraycolsep=1pt
\begin{array}{lll}
 (-\Delta)^{\alpha} u(x)=f(u(x))+g(x),\ \ \ \ &
x\in \Omega,\\[2mm]
 u(x)=0,& x\in \Omega^c.
\end{array}
\right.
\end{equation}
Then $u$ is symmetric with respect to $x_1$ and it is strictly decreasing in the $x_1$ direction
for $x=(x_1,x')\in\Omega$, $x_1>0$.
\end{teo}

\setcounter{equation}{0}
\section{Symmetry of solutions  in $\R^N$}

In this section we study radial symmetry results for positive
solution of equation \equ{eq wyq1} in $\R^N$, in particular we will
provide a proof of Theorem \ref{teo r1}. In the case of the whole
space, the moving planes procedure needs to be started in a
different way, because we cannot use the Maximum Principle for small
domains. We use the moving plane method as for the second order
equation as in the work by Li \cite{LCM} (see also \cite{PR}).

In this section we use the notation introduced  in \equ{d1}-\equ{d3}
and we let $u$ be a classical positive solution of  \equ{eq wyq1}.
In order to prove Theorem \ref{teo r1} we need some preliminary
lemmas.
\begin{lemma}\label{lemma sev1}
Under the assumptions of Theorem \ref{teo r1}, for any  $\lambda\in\R$, we have
$$
\int_{\Sigma_\lambda} (f(u_\lambda)-f(u))^+(u_\lambda-u)^+dx
<+\infty.
$$
\end{lemma}
\noindent{\bf Proof.}
By our hypothesis, for any given $\lambda\in\R$, we may choose  $R>1$  and some constant $c>1$ such that
$$ \frac{1}{c|x|^m}\leq u(x), u_\lambda(x)\leq\frac{c}{|x|^m}<s_0 \ \ \ for \ all\ x\in B^c_R,$$
where $s_0$ is the constant in condition (F2).

If $u_\lambda(x)> u(x)$ for some $x\in\Sigma_\lambda\cap B^c_R,$  we have
$0<u(x)<u_\lambda(x)<s_0$. Using (\ref{eq rq2}) with $v=u_\lambda(x)$, then
$$ \frac{f(u_\lambda(x))-f(u(x))}{u_\lambda(x)-u(x)}
\leq C(u(x)+u_\lambda(x))^\gamma\leq 2^\gamma Cu^\gamma_\lambda(x),$$
then
\begin{eqnarray*}
(f(u_\lambda(x))-f(u(x)))^+(u_\lambda(x)-u(x))^+
&\leq& 2^\gamma Cu^\gamma_\lambda(x)[(u_\lambda(x)-u(x))^+]^2
\\&\leq& \tilde{C}u^{\gamma+2}_\lambda(x),
\end{eqnarray*}
for certain $\tilde{C}>0$. We observe that, if $u_\lambda(x)\leq u(x)$ for some $x\in\Sigma_\lambda\cap B^c_R,$
then inequality above is obvious. Therefore,
\begin{eqnarray*}
(f(u_\lambda)-f(u))^+(u_\lambda-u)^+\leq \tilde{C}u^{\gamma+2}_\lambda \ \ \  in\ \ \Sigma_\lambda\cap B^c_R.
 \end{eqnarray*}
Now we integrate in $\Sigma_\lambda\cap B^c_R$ to obtain
\begin{eqnarray*}
\int_{\Sigma_\lambda\cap B^c_R}(f(u_\lambda)-f(u))^+(u_\lambda-u)^+dx
&\leq&  \tilde{C}\int_{\Sigma_\lambda\cap {B^c_R}}u^{\gamma+2}_\lambda(x)dx
\\&\leq&  C\int_{ {B^c_R}}|x|^{-m(\gamma+2)}dx<+\infty,
\end{eqnarray*}
where the last inequality holds by \equ{m}. Since $u$ and $u_\lambda$ are bounded and $f$ is locally Lipschitz, we have
\begin{eqnarray*}
 \int_{\Sigma_\lambda\cap B_R}(f(u_\lambda)-f(u))^+(u_\lambda-u)^+dx
<+\infty
\end{eqnarray*}
and the proof is complete.
 \hfill$\Box$\\

It will be convenient for our analysis to define the following function
\begin{equation}\label{eq dem5}
w(x)=\left\{ \arraycolsep=1pt
\begin{array}{lll}
(u_\lambda-u)^+(x),\ \ \ \ &
x\in \Sigma_\lambda,\\[2mm]
(u_\lambda-u)^-(x),\ \ \ \ &
x\in \Sigma_\lambda^c,
\end{array}
\right.
\end{equation}
where $(u_\lambda-u)^+(x)=\max\{(u_\lambda-u)(x),\ 0\}$, $(u_\lambda-u)^-(x)=\min\{(u_\lambda-u)(x),\ 0\}$.
We have
\begin{lemma}\label{lemma dem4}
Under the assumptions of Theorem  \ref{teo r1}, there exists a constant $C>0$ such that
\begin{equation}\label{eq dem6}
\int_{\Sigma_\lambda} (-\Delta)^{\alpha}({u_\lambda}-u) (u_\lambda-u)^+dx
\geq C (\int_{\Sigma_\lambda}|w|^{\frac{2N}{N-2\alpha}}dx)^{\frac{N-2\alpha}{N}}.
\end{equation}
\end{lemma}
\noindent{\bf Proof.}
 We start observing that, given $x\in\Sigma_\lambda$, we have
\begin{eqnarray*}
w(x_\lambda)&=&(u_\lambda-u)^-(x_\lambda)=\min\{(u_\lambda-u)(x_\lambda),\ 0\}=\min\{(u-u_\lambda)(x),\ 0\}
\\&=&-\max\{(u_\lambda-u)(x),\ 0\}=-(u_\lambda-u)^+(x)
= -w(x)
\end{eqnarray*}
and similarly
 $w(x)=-w(x_\lambda)$ for $x\in\Sigma_\lambda^c$ so that
 \begin{equation}\label{signo}
 w(x)=-w(x_\lambda)\quad\mbox{for}\quad  x\in \R^N.
 \end{equation}
This implies
 \begin{eqnarray}\label{doble}
\int_{\R^N}|w|^{\frac{2N}{N-2\alpha}}dx=\int_{\Sigma_\lambda}|w|^{\frac{2N}{N-2\alpha}}dx+\int_{\Sigma_\lambda^c}|w|^{\frac{2N}{N-2\alpha}}dx
=2\int_{\Sigma_\lambda}|w|^{\frac{2N}{N-2\alpha}}dx.
\end{eqnarray}
Next we see that for any  $x\in \Sigma_\lambda\cap {\rm{supp}}(w)$ we have  that $w(x)=(u_\lambda-u)(x)$  and
$$(-\Delta)^{\alpha}(u_\lambda-u)(x)\ge (-\Delta)^{\alpha}w(x),\quad \forall\ x\in \Sigma_\lambda\cap \rm{supp}(w), $$
\begin{eqnarray}
&&(-\Delta)^{\alpha}w(x)- (-\Delta)^{\alpha}(u_\lambda-u)(x)
= \int_{\R^N}\frac{(u_\lambda-u)(z)-w(z)}{|x-z|^{N+2\alpha}}dz  \nonumber
\\&=&\int_{\Sigma_\lambda\cap({\rm{supp}}(w))^c}\frac{(u_\lambda-u)(z)}{|x-z|^{N+2\alpha}}dz
+\int_{\Sigma_\lambda^c\cap({\rm{supp}}(w))^c}\frac{(u_\lambda-u)(z)}{|x-z|^{N+2\alpha}}dz\nonumber
\\&=&\int_{\Sigma_\lambda\cap({\rm{supp}}(w))^c}(u_\lambda-u)(z)
(\frac1{|x-z|^{N+2\alpha}}-\frac1{|x-z_\lambda|^{N+2\alpha}})dz
\le 0, \label{desiw}
\end{eqnarray}
where we used that  $u_\lambda-u\leq0$ in $\Sigma_\lambda\cap({\rm{supp}}(w))^c$ and
$|x-z|\leq |x-z_\lambda|$ for $x,z\in\Sigma_\lambda.$
From \equ{desiw}, using the equation and Lemma \ref{lemma sev1} we find that
\begin{eqnarray}
 \int_{\Sigma_\lambda} (-\Delta)^{{\alpha}}w\, wdx &\le&  \int_{\Sigma_\lambda} (-\Delta)^{{\alpha}}(u_\lambda-u)(u_\lambda -u)^+dx\label{cotaD0}\\&\le& \int_{\Sigma_\lambda} (f(u_\lambda)-f(u))^+ (u_\lambda -u)^+dx<\infty.
 \label{cotaD}
  \end{eqnarray}
From here the following integrals are finite and,  taking into account \equ{signo},  we obtain that
 \begin{eqnarray}
\int_{\R^N}|(-\Delta)^{\frac{\alpha}{2}}w|^2dx&=&
\int_{\Sigma_\lambda}|(-\Delta)^{\frac{\alpha}{2}}w|^2dx+\int_{\Sigma_\lambda^c}|(-\Delta)^{\frac{\alpha}{2}}w|^2dx
\nonumber \\&=&2\int_{\Sigma_\lambda}|(-\Delta)^{\frac{\alpha}{2}}w|^2dx. \label{cotaE}
\end{eqnarray}
Now we can use the Sobolev embedding from $H^\alpha(\R^N)$ to
$L^{\frac{2N}{N-2\alpha}}(\R^N)$ to find  a constant $C$ so that
 \begin{eqnarray}
 \int_{\Sigma_\lambda}|(-\Delta)^{\frac{\alpha}{2}}w|^2dx
&=&\frac12 \int_{\R^N}|(-\Delta)^{\frac{\alpha}{2}}w|^2dx
\geq C (\int_{\R^N}|w|^{\frac{2N}{N-2\alpha}}dx)^{\frac{N-2\alpha}{N}}\nonumber \\
&=&C(2\int_{\Sigma_\lambda}|w|^{\frac{2N}{N-2\alpha}}dx)^{\frac{N-2\alpha}{N}}.
\label{sobolev}\end{eqnarray}
On the other hand, from \equ{signo} and \equ{cotaD0} we find that
\begin{eqnarray}
 \int_{\R^N}|(-\Delta)^{\frac{\alpha}{2}}w|^2dx \nonumber
 &=&\int_{\R^N}(-\Delta)^{\alpha}w\cdot wdx
 =2\int_{\Sigma_\lambda}(-\Delta)^{\alpha}w\cdot wdx
 \\&\leq&2 \int_{\Sigma_\lambda} (-\Delta)^{\alpha}({u_\lambda}-u) (u_\lambda-u)^+dx. \label{otra}
 \end{eqnarray}
From \equ{sobolev} and \equ{otra}  the proof of the lemma is completed.
%
 \hfill$\Box$\\

\medskip

Now we are ready to complete the

\noindent{\bf Proof of Theorem \ref{teo r1}.}
 We divide the proof into three steps.

 \noindent\textbf{Step 1:} We  show that $\lambda_0:=\sup\{\lambda\ |\ u_\lambda\leq u\ in\ \Sigma_\lambda \}$ is finite. Using $(u_\lambda-u)^+$ as a test function in the equation for
$u$ and $u_\lambda$,  using
(\ref{eq rq2}) and  H\"{o}lder inequality, for $\lambda$ big (negative), we find that
\begin{eqnarray*}
\int_{\Sigma_\lambda} (-\Delta)^{\alpha}({u_\lambda}-u) (u_\lambda-u)^+dx
\!\!\! \!\!\! \!\!\!&&= \int_{\Sigma_\lambda} (f(u_\lambda)-f(u))(u_\lambda-u)^+dx
\\&&\leq \int_{\Sigma_\lambda}[\frac{f(u_\lambda)-f(u)}{u_\lambda-u}]^+[(u_\lambda-u)^+]^2dx
\\&&\leq C\int_{\Sigma_\lambda}{u^\gamma_\lambda}w^2dx
\leq \bar{C}\int_{\Sigma_\lambda}|x_\lambda|^{-m\gamma}w^2dx
\\&&\leq \bar{C}(\int_{\Sigma_\lambda}|x_\lambda|^{-\frac{Nm\gamma}{2\alpha}}dx)^{\frac{2\alpha}{N}}
  (\int_{\Sigma_\lambda}|w|^{\frac{2N}{N-2\alpha}}dx)^{\frac{N-2\alpha}{N}}.
\end{eqnarray*}
By  Lemma \ref{lemma dem4}, there exists a constant $C>0$ such that
\begin{eqnarray*}
(\int_{\Sigma_\lambda}|w|^{\frac{2N}{N-2\alpha}}dx)^{\frac{N-2\alpha}{N}}
&\leq& C(\int_{\Sigma_\lambda}|x_\lambda|^{-\frac{Nm\gamma}{2\alpha}}dx)^{\frac{2\alpha}{N}}
  (\int_{\Sigma_\lambda}|w|^{\frac{2N}{N-2\alpha}}dx)^{\frac{N-2\alpha}{N}},
\end{eqnarray*}
but we have
\begin{eqnarray*}
\int_{\Sigma_\lambda}|x_\lambda|^{-\frac{Nm\gamma}{2\alpha}}dx&\le& \int_{\Sigma^c_\lambda}|x|^{-\frac{Nm\gamma}{2\alpha}}dx
\leq \int_{B^c_{|\lambda|}}|x|^{-\frac{Nm\gamma}{2\alpha}}dx
= c {|\lambda|}^{\frac{N}{2\alpha}(2\alpha-m\gamma)},
\end{eqnarray*}
so that, using (\ref{m}), we can choose  $R>0$ big enough such that  $ C R^{2\alpha-m\gamma}\leq\frac{1}{2}$,
then we obtain
$$\int_{\Sigma_\lambda}|w|^{\frac{2N}{N-2\alpha}}dx=0,\ \  \forall \ \lambda<-R.$$
Thus $w=0$ in $\Sigma_\lambda$ and then $u_\lambda\leq u$  in $\Sigma_\lambda,$ for all $\lambda<-R,$ concluding that $\lambda_0\geq -R.$
On the other hand, since $u$ decays at infinity, then there exists $\lambda_1$ such that
$u(x)<u_{\lambda_1}(x)$ for some $x\in \Sigma_{\lambda_1}.$
 Hence  $\lambda_0$ is finite.
 \medskip

\noindent\textbf{Step 2:}
We prove that $u\equiv u_{\lambda_0}$ in $\Sigma_{\lambda_0}$.
Assuming the contrary, we have  $u\neq u_{\lambda_0}$ and  $u\geq u_{\lambda_0}$ in $\Sigma_{\lambda_0}$.
Assume next that there exists $x_0\in \Sigma_{\lambda_0}$ such that
$u_{\lambda_0}(x_0)=u(x_0),$ then
 we have
\begin{equation}\label{eq abcd}
(-\Delta)^{\alpha} u_{\lambda_0}(x_0)-(-\Delta)^{\alpha} u(x_0)
=f(u_{\lambda_0}(x_0))-f(u(x_0))=0.
\end{equation}
On the other hand,
\begin{eqnarray*}
&&(-\Delta)^{\alpha} u_{\lambda_0}(x_0)-(-\Delta)^{\alpha} u(x_0)
=-\int_{\R^N}\frac{u_{\lambda_0}(y)-u(y)}{|x_0-y|^{N+2\alpha}}dy
\\&=&-\int_{\Sigma_{\lambda_0}}(u_{\lambda_0}(y)-u(y))(\frac{1}{|x_0-y|^{N+2\alpha}}-\frac{1}{|x_0-y_{\lambda_0}|^{N+2\alpha}})dy
>0,
\end{eqnarray*}
which contradicts  (\ref{eq abcd}). As a sequence, $u> u_{\lambda_0}$ in $\Sigma_{\lambda_0}$.

To complete Step 2, we only need to prove that
$u\geq u_\lambda$ in $\Sigma_\lambda$ continues to hold when ${\lambda_0}<\lambda<{\lambda_0}+\varepsilon$,
where $\varepsilon>0$ small.
Let us consider then $\varepsilon>0$, to be chosen later, and
take  $\lambda\in({\lambda_0},{\lambda_0}+\varepsilon)$. Let $P=(\lambda,0)$ and $B(P,R)$ be the ball centered at $P$ and with radius $R>1$ to be chosen later.
 Define $\tilde{B}=\Sigma_\lambda\cap B(P,R)$ and let us consider  $(u_\lambda-u)^+$  test function in the equation for
$u$ and $u_\lambda$ in $\Sigma_\lambda$, then from Lemma \ref{lemma dem4} we find
\begin{eqnarray}\label{E0}
(\int_{\Sigma_\lambda}|w|^{\frac{2N}{N-2\alpha}}dx)^{\frac{N-2\alpha}{N}}
&\le& C\int_{\Sigma_\lambda } (f(u_\lambda)-f(u))(u_\lambda-u)^+dx.
\end{eqnarray}
We estimate the integral on the right.
Since $f$ is locally Lipschitz, using H\"older inequality, we have
\begin{eqnarray}
\nonumber && \int_{\tilde{B}} (f(u_\lambda)-f(u))(u_\lambda-u)^+dx
\leq C\int_{\tilde{B}}|w|^2\chi_{{\rm{supp}}{(u_\lambda-u)^+}}dx
\\&=& C|\tilde{B}\cap{{\rm{supp}}{(u_\lambda-u)^+}}|^{\frac{2\alpha}{N}}
  (\int_{\tilde{B}}|w|^{\frac{2N}{N-2\alpha}}dx)^{\frac{N-2\alpha}{N}}.\label{E1}
\end{eqnarray}
On the other hand, for the integral over  $\Sigma_\lambda\setminus{\tilde{B}}$, we assume $R$ and $R_0$ are such that
 $\Sigma_\lambda\setminus{\tilde{B}}\subset {B^c(P,R)}\subset B^c_{R_0}(0)$,
proceeding as in Step 1,  we have
\begin{eqnarray}
\int_{\Sigma_\lambda\setminus{\tilde{B}}} (f(u_\lambda)-f(u))(u_\lambda-u)^+dx
&\leq& {C}\int_{\Sigma_\lambda\setminus{\tilde{B}}}u^\gamma_\lambda w^2dx
\nonumber \\&\leq& {C}(\int_{\Sigma_\lambda\setminus{\tilde{B}}}|x_\lambda|^{-\frac{Nm\gamma}{2\alpha}}dx)^{\frac{2\alpha}{N}}
  (\int_{\Sigma_\lambda}|w|^{\frac{2N}{N-2\alpha}}dx)^{\frac{N-2\alpha}{N}}\nonumber
\\&\le& C {R_0}^{2\alpha-m\gamma} (\int_{\Sigma_\lambda}|w|^{\frac{2N}{N-2\alpha}}dx)^{\frac{N-2\alpha}{N}}.
\label{E2}\end{eqnarray}
Now we choose $R_0$ such that $C {R_0}^{2\alpha-m\gamma}<1/2$, then choose $R$ so that $\Sigma_\lambda\setminus{\tilde{B}}\subset {B^c(P,R)}\subset B^c_{R_0}(0)$ and then choose
$\varepsilon>0$ so that
$C|\tilde{B}\cap{{\rm{supp}}{(u_\lambda-u)^+}}|^{\frac{2\alpha}{N}}<1/2$. With this choice of the parameters, from \equ{E0}, \equ{E1} and \equ{E2} it follows that
%
%
 $w=0$ in $\Sigma_\lambda$, which is a contradiction, completeing Step 2.
 \medskip

\noindent\textbf{Step 3:}
By translation, we may say that $\lambda_0=0.$ An repeating the argument from the other side, we find that $u$ is symmetric about $x_1$-axis. Using
the same argument in any arbitrary direction, we finally conclude that $u$ is radially symmetric.

Finally, we prove that $u(r)$ is strictly decreasing in  $r>0$, by using the same arguments as in the case of a ball.
This completes the proof.
 \hfill$\Box$\\

At the end of this section we want to give a theorem on radial symmetry of solutions for equation  \equ{eq wyq1}  in a case where $f$ is only locally Lipschitz in $(0,\infty)$, see \cite{CEF2} and \cite{CEF} for  the case of the Laplacian.
In precise terms we have
\begin{teo}\label{teo 42}
Let $u$ be a positive classical solution of
\begin{equation}\label{eq add42}
\left\{ \arraycolsep=1pt
\begin{array}{lll}
 (-\Delta)^{\alpha} u=u^p-u^q\ \ \ in\ \ \R^N,\\[2mm]
 u>0\ \ in\ \ \R^N,\ \ \lim_{|x|\to\infty}u(x)=0,
\end{array}
\right.
\end{equation}
satisfying
\begin{equation}\label{eq add43}
 u(x)=O(|x|^{-\frac{N+2\alpha}{q}})\ \ \ \ as\  |x|\to\infty,
 \end{equation}
where $\alpha\in(0,1),$ $N\geq2$ and $0<q<1<p$.
 Then $u$ is radially symmetric
and strictly decreasing about some point.
\end{teo}
\noindent{\bf Proof.}
We denote $f(u)=u^p-u^q$ for $u>0$, and consider  $\gamma>0$ and $s_0$ small enough,
then for all $u,v$ satisfying $0<u<v<s_0$, we have
$$\frac{f(v)-f(u)}{v-u}
<0\leq C(u+v)^\gamma,$$ for some constant $C>0$, so that (F2) holds.
We also observe that for a positive classical solution $u$ of
(\ref{eq add42}), $u\geq c$ in any bounded domain $\Omega$, for a
constant $c>0$ depending on $\Omega$ and then, in (\ref{E1}) we may
use Lipschitz continuity of $f$ in the bounded interval $[c,\sup
u]$. We set $m=\frac{N+2\alpha}{q}$ and $\gamma$ may be chosen so
that \equ{m} holds. The  proof of Theorem \ref{teo 42} goes in the
same way as that of Theorem \ref{teo r1}.
 \hfill$\Box$

\begin{remark}
In a work by
Valdebenito \cite{valdebenito}, the estimate \equ{eq add43} is obtained by using super solutions and  Theorem  \ref{teo 42}  is proved using the local extension of equation  \equ{eq add42} as given by Caffarelli and Silvestre in \cite{CS} and then using a regular moving planes argument as developed for elliptic equations with non-linear boundary conditions  by Terracini \cite{terracini}.
\end{remark}

\setcounter{equation}{0}
\section{Symmetry results for system}

The aim of this section is to prove Theorem \ref{teo 2} by the
moving planes method applied to a system of equations in the unit
ball $B_1$. Let $\Sigma_\lambda$ and $T_\lambda$ be defined as in
Section \S 3. For $x=(x_1,x')\in\R^N$ and $\lambda\in (0,1)$ we let
$x_\lambda=(2\lambda-x_1,x')$,
$$u_\lambda(x)=u(x_\lambda),\ \ \  \ w_{\lambda,u}(x)=u_\lambda(x)-u(x),$$
$$v_\lambda(x)=v(x_\lambda),\quad\mbox{and}\quad   w_{\lambda,v}(x)=v_\lambda(x)-v(x).$$

\medskip

\noindent{\bf Proof of  Theorem \ref{teo 2}.}
 We will split this proof into three steps.

\noindent\textbf{Step 1:}
We start the moving planes proving that  if $\lambda$ is close to $1$, then $w_{\lambda,u}$ and $w_{\lambda,v}$ are positive in $\Sigma_\lambda$. For that purpose we define
$$
\Sigma_{\lambda,u}^-=\{x\in \Sigma_\lambda\ |\ w_{\lambda,u}(x)<0\}\quad\mbox{and}\quad
\Sigma_{\lambda,v}^-=\{x\in \Sigma_\lambda\ |\ w_{\lambda,v}(x)<0\}.$$
We show next that  $\Sigma_{\lambda,u}^-$ is empty for $\lambda$ close to 1.
Assume, by contradiction, that $\Sigma_{\lambda,u}^-$ is not empty and define
\begin{equation}\label{eq 25}
w_{\lambda,u}^+(x)=\left\{ \arraycolsep=1pt
\begin{array}{lll}
w_{\lambda,u}(x),\ \ \ \ &
x\in \Sigma_{\lambda,u}^-,\\[2mm]
0,& x\in \R^N\setminus \Sigma_{\lambda,u}^-
\end{array}
\right.
\end{equation}
and
\begin{equation}\label{eq 26}
w_{\lambda,u}^-(x)=\left\{ \arraycolsep=1pt
\begin{array}{lll}
0,\ \ \ \ &
x\in \Sigma_{\lambda,u}^-,\\[2mm]
w_{\lambda,u}(x),\ \ \ \ & x\in \R^N\setminus \Sigma_{\lambda,u}^-.
\end{array}
\right.
\end{equation}
Using the arguments given in Step 1 of the proof of Theorem \ref{teo 1}, we get
\begin{equation}\label{eq 270}
(-\Delta)^{\alpha_1} w_{\lambda,u}^+(x)\ge (-\Delta)^{\alpha_1} w_{\lambda,u}(x)
\quad\mbox{and}\quad (-\Delta)^{\alpha_1} w_{\lambda,u}^-(x)\le0,
\end{equation}
for all $x\in\Sigma^-_{\lambda,u}$.
From here, using equation  (\ref{eq 21}), for $x\in\Sigma^-_{\lambda,u}$ we have
\begin{eqnarray}
(-\Delta)^{\alpha_1} w_{\lambda,u}^+(x)
&\ge& (-\Delta)^{\alpha_1} u_\lambda(x)-(-\Delta)^{\alpha_1} u(x)
\nonumber \\&=& f_1(v_\lambda(x))+g_1(x_\lambda)-f_1(v(x))-g_1(x)
\nonumber\\&=&\varphi_v(x)
w_{\lambda,v}(x)+g_1(x_\lambda)-g_1(x)
\nonumber \\&\ge &\varphi_v(x)
w_{\lambda,v}(x), \label{dess}
\end{eqnarray}
where
 $\varphi_v(x)=({f_1(v_\lambda(x))-f_1(v(x)))}/({v_\lambda(x)-v(x)})$ and  where we used that $g_1$ is radially symmetric and decreasing, with $|x|>|x_\lambda|$.
 We further observe that, since $f_1$ is locally Lipschitz continuous, we have that $\varphi_v(\cdot)\in
L^\infty(\Sigma^-_{\lambda,u})$.
Now we consider  \equ{dess} together with $w_{\lambda,u}^+=0$ in $(\Sigma_{\lambda,u}^-)^c$ and   $w_{\lambda,u}^+<0$ in  $\Sigma_{\lambda,u}^-$, to  use
 Proposition \ref{pro abp1} to find a
constant $C>0$, depending on $N$ and $\alpha$ only, such that
\begin{equation}\label{eq g5}
\|w_{\lambda,u}^+\|_{L^\infty(\Sigma_{\lambda,u}^-)} \le
C
\|(-\varphi_v
w_{\lambda,v})^+\|^{1-\alpha_1}_{L^\infty(\Sigma_{\lambda,u}^-)}
\|(-\varphi_v
w_{\lambda,v})^+\|^{\alpha_1}_{L^N(\Sigma_{\lambda,u}^-)}
 \end{equation}
We observe that  $diam (\Sigma_{\lambda,u}^-)\le 1.$
 Since $f_1$ is increasing, we have
 \begin{eqnarray}\label{eq g4}
-\varphi_v w_{\lambda,v}&=&f_1(v)-f_1(v_\lambda)\leq0\ \ in\
(\Sigma_{\lambda,v}^-)^c\quad \mbox{and}\\
-\varphi_v
w_{\lambda,v}&=&f_1(v)-f_1(v_\lambda)>0\ \ in\ \Sigma_{\lambda,v}^-.\label{eq g41}
\end{eqnarray}
Denoting $\Sigma_{\lambda}^-=\Sigma_{\lambda,u}^-\cap \Sigma_{\lambda,v}^-,$ from (\ref{eq g5}), (\ref{eq g4}) and (\ref{eq g41}),  we obtain
\begin{eqnarray}\label{iv}
\|w_{\lambda,u}^+\|_{L^\infty(\Sigma_{\lambda,u}^-)} \le
C
\|(-\varphi_v
w_{\lambda,v})^+\|_{L^\infty(\Sigma_{\lambda}^-)}
|\Sigma_{\lambda}^-|^{\frac{\alpha_1}{N}},
\end{eqnarray}
 Similar to \equ{eq 25} and  \equ{eq 26}, we define
\begin{eqnarray*}
w_{\lambda,v}^+(x)=\left\{ \arraycolsep=1pt
\begin{array}{lll}
w_{\lambda,v}(x),\ \ \ \ &
x\in \Sigma_{\lambda,v}^-,\\[2mm]
0,& x\in \R^N\setminus \Sigma_{\lambda,v}^-
\end{array}
\right.
\end{eqnarray*}
and
\begin{eqnarray*}
w_{\lambda,v}^-(x)=\left\{ \arraycolsep=1pt
\begin{array}{lll}
0,\ \ \ \ &
x\in \Sigma_{\lambda,v}^-,\\[2mm]
w_{\lambda,v}(x),\ \ \ \ & x\in \R^N\setminus \Sigma_{\lambda,v}^-.
\end{array}
\right.
\end{eqnarray*}
With this definition \equ{iv} becomes
\begin{equation}\label{eq g1}
\|w_{\lambda,u}^+\|_{L^\infty(\Sigma_{\lambda,u}^-)} \le
C
\|w_{\lambda,v}^+\|
_{L^\infty(\Sigma_{\lambda}^-)}|\Sigma_{\lambda}^-|^{\frac{\alpha_1}{N}},
\end{equation}
where we used that $\varphi_v$ is bounded and we have changed the constant $C$, if necessary. At this point we observe that if $w_{\lambda,v}^+=0$  then $w_{\lambda,u}^+=0$ providing a contradiction. Thus we have that
$\Sigma_{\lambda,v}^-\not=\emptyset$
and we may argue in a completely analogous way to obtain
\begin{equation}\label{eq g2}
\|w_{\lambda,v}^+\|_{L^\infty(\Sigma_{\lambda,v}^-)}\le
C
\|w_{\lambda,u}^+\|
_{L^\infty(\Sigma_{\lambda}^-)}|\Sigma_{\lambda}^-|^{\frac{\alpha_2}{N}},
\end{equation}
that combined with  \equ{eq g1} yields
%
$$\|w_{\lambda,u}^+\|_{L^\infty(\Sigma_{\lambda,u}^-)}
\le C^2
|\Sigma_{\lambda}^-|^\frac{{\alpha_1}+{\alpha_2}}{N}\|w_{\lambda,u}^+\|_{L^\infty(\Sigma_{\lambda,u}^-)},
$$
and
$$\|w_{\lambda,v}^+\|_{L^\infty(\Sigma_{\lambda,v}^-)}
\le C^2
|\Sigma_{\lambda}^-|^\frac{{\alpha_1}+{\alpha_2}}{N}\|w_{\lambda,v}^+\|_{L^\infty(\Sigma_{\lambda,v}^-)}.
$$
Now we just take
$\lambda$ close enough to $1$ so that $C^2|\Sigma_{\lambda}^-|^\frac{{\alpha_1}+{\alpha_2}}{N} <1$ and we conclude that
$\|w_{\lambda,u}^+\|_{L^\infty(\Sigma_{\lambda,u}^-)}=\|w_{\lambda,v}^+\|_{L^\infty(\Sigma_{\lambda,v}^-)}=0,$
so $|\Sigma_{\lambda,u}^-|=|\Sigma_{\lambda,v}^-|=0$ and since
$\Sigma_{\lambda,u}^-$ and $\Sigma_{\lambda,v}^-$  are open we have that
$\Sigma_{\lambda,u}^-,\Sigma_{\lambda,v}^-=\O$, which is a contradiction.

Thus we have that
$w_{\lambda,u}\ge0$ in $\Sigma_{\lambda}$ when $\lambda$ is close enough to
$1$. Similarly, we obtain $w_{\lambda,v}\ge0$ in $\Sigma_{\lambda}$
for $\lambda$ close to $1$. In order to complete Step 1 we will prove  a bit more general statement that will be useful later, that is,
given  $0<\lambda<1$, if $w_{\lambda,u}\ge 0, w_{\lambda,v}\geq0$, $w_{\lambda,u}\not\equiv0$ and  $w_{\lambda,v}\not\equiv0$ in $\Sigma_\lambda$,
then  $w_{\lambda,u}>0$ and $w_{\lambda,v}>0$ in $\Sigma_\lambda$. For proving this property suppose  there exists $x_0\in \Sigma_\lambda$ such that
\begin{equation}\label{eq f1}
w_{\lambda,u}(x_0)=0.
\end{equation}
On one hand, by using similar arguments yielding \equ{ff} we find that
\begin{equation}\label{eq f2}
(-\Delta)^{\alpha_1} w_{\lambda,u}(x_0)<0.
\end{equation}
On the other hand, by our assumption we have that
 $ w_{\lambda,v}(x_0)=v_\lambda(x_0)-v(x_0)\geq0$ and since $|x_0|>|(x_0)_\lambda|$, from the monotonicity hypothesis on  $f_1$ and  $g_1$, we obtain
 $$f_1(v_\lambda(x_0))\geq f_1(v(x_0)),\ \ \ \ \  g_1((x_0)_\lambda)\geq g_1(x_0).$$
 Thus, using  (\ref{eq 21}), we find
\begin{eqnarray*}
(-\Delta)^{\alpha_1} w_{\lambda,u}(x_0)
&=&
f_1(v_\lambda(x_0))+g_1((x_0)_\lambda)-f_1(v(x_0))-g_1(x_0)\ge 0,
\end{eqnarray*}
which is impossible with (\ref{eq f2}). This completes Step 1.

\noindent\textbf{Step 2:}
 We prove that $\lambda_0=0$, where
$$\lambda_0=\inf\{\lambda\in(0,1)\ |\  w_{\lambda,u}\ ,\ w_{\lambda,v}>0\ \ \rm{in}\ \ \Sigma_\lambda\}.$$
If not, that is, if $\lambda_0>0$ we
have that $w_{\lambda_0,u},w_{\lambda_0,v}\ge0$ and
$w_{\lambda_0,u},w_{\lambda_0,v}\not\equiv0$ in $\Sigma_{\lambda_0}$. If we use the property we just proved above, we may assume that
 $w_{\lambda_0,u}>0$ and $w_{\lambda_0,v}>0$ in $\Sigma_{\lambda_0}$. In what follows we argue that
 the plane can be moved to left, that is, that   there exists $\epsilon\in(0,\lambda)$ such that
$w_{{\lambda_\epsilon},u}>0$ and $w_{{\lambda_\epsilon},v}>0$ in $\Sigma_{\lambda_\epsilon}$,
where $\lambda_\epsilon=\lambda_0-\epsilon$, providing a contradiction with the definition of $\lambda_0$.

Let us consider the set $D_\mu=\{x\in\Sigma_\lambda\ | \ dist(x,\partial\Sigma_\lambda)\ge \mu\}$ for $\mu>0$ small.
Since $w_{\lambda,u},w_{\lambda,v}>0$ in $\Sigma_\lambda$ and $D_\mu$ is compact,
then there exists $\mu_0>0$ such that $w_{\lambda,u},w_{\lambda,v}\ge \mu_0$ in $D_\mu$.
By continuity of $w_{\lambda,u}(x)$ and $w_{\lambda,v}(x)$, for $\epsilon>0$ small enough,
 we have that
$$w_{\lambda_\epsilon,u},\ w_{\lambda_\epsilon,v}\ge0\ \ \rm{in}\ \ D_\mu$$
and, as a consequence,
$\Sigma_{\lambda_\epsilon,u}^-, \Sigma_{\lambda_\epsilon,v}^- \subset\Sigma_{\lambda_\epsilon}\setminus D_\mu,$
and $|\Sigma_{\lambda_\epsilon,u}^-|$ and $|\Sigma_{\lambda_\epsilon,v}^-|$ are small if $\epsilon$ and $\mu$ are small.

Since $f_1$ and $f_2$ are locally Lipschitz continuous and increasing, $g_1$ and $g_2$ are
radially symmetric and decreasing, we may repeat the arguments given in Step 1 to obtain
%
$$\|w_{\lambda_\epsilon,u}^+\|_{L^\infty(\Sigma_{\lambda_\epsilon,u}^-)}\le
C^2
|\Sigma_{\lambda_\epsilon}^-|^\frac{{\alpha_1}+{\alpha_2}}{N}\|w_{\lambda_\epsilon,u}^+\|_{L^\infty(\Sigma_{\lambda_\epsilon,u}^-)}
$$
and
$$\|w_{\lambda_\epsilon,v}^+\|_{L^\infty(\Sigma_{\lambda_\epsilon,v}^-)}\le
C^2
|\Sigma_{\lambda_\epsilon}^-|^\frac{{\alpha_1}+{\alpha_2}}{N}\|w_{\lambda_\epsilon,v}^+\|_{L^\infty(\Sigma_{\lambda_\epsilon,v}^-)}
$$
where
$\Sigma_{\lambda_\epsilon}^-=\Sigma_{\lambda_\epsilon,u}^-\cap
\Sigma_{\lambda_\epsilon,v}^-$.
Now we may choose
  $\epsilon$ and $\mu$
 small such that
$C^2
|\Sigma_{\lambda_\epsilon}^-|^\frac{{\alpha_1}+{\alpha_2}}{N}<1,$
then we obtain
$\|w_{\lambda_\epsilon,u}^+\|_{L^\infty(\Sigma_{\lambda_\epsilon,u}^-)}=\|w_{\lambda_\epsilon,v}^+\|_{L^\infty(\Sigma_{\lambda_\epsilon,v}^-)}=0$. From here we argue as in Step 1 to obtain that
$w_{\lambda_\epsilon,u}$ and $w_{\lambda_\epsilon,v}$ are positive
in $\Sigma_{\lambda_\epsilon}$, completing Step 2.

Finally, we obtain that $u$ and $v$ are
radially symmetric and strictly decreasing respect to  $r=|x|$ for
$r\in(0,1)$ in the same way in Step 3 in the proof of Theorem \ref{teo 1}.
 \hfill$\Box$\\

\setcounter{equation}{0}
\section{The case of a non-local operator  with non-homogeneous kernel.}

The main purpose of this section is to discuss radial symmetry for a problem with a non-local    operator $\mathcal{L}$ of fractional order, but with a non-homogeneous kernel. The operator is  defined as follows:
\begin{equation}\label{eq 2}
\mathcal{L} u(x)=P.V.\int_{\R^N}(u(x)-u(y)){K_\mu}(x-y)dy,
\end{equation}
where the kernel ${K_\mu}$ satisfies that
\begin{equation}\label{eq 120}
K_\mu(x)=\left\{ \arraycolsep=1pt
\begin{array}{lll}
 \frac1{|x|^{N+2\alpha_1}},\ \ \ \ & |x|<1,\\[2mm]
 \frac{\mu}{|x|^{N+2\alpha_2}},\ \ \ \ & |x|\geq 1
\end{array}
\right.
\end{equation}
with $\mu\in[0,1]$ and $\alpha_1,\alpha_2\in(0,1)$.
Being more precise, we consider the equation
\begin{equation}\label{eq 51}
\left\{ \arraycolsep=1pt
\begin{array}{lll}
\mathcal{L} u(x)=f(u(x))+g(x),\ \ \ \ &
x\in B_1,\\[2mm]
 u(x)=0,& x\in B_1^c,
\end{array}
\right.
\end{equation}
and our theorem states
\begin{teo}\label{teo 5}
 Assume that the
function  $f$ satisfies $(F1)$ and $g$ satisfies $(G)$. If $u$ is a positive classical  solution of (\ref{eq 51}),
then $u$ must be radially symmetric and strictly decreasing in
$r=|x|$ for $r\in(0,1)$.
\end{teo}

The idea for Theorem \ref{teo 5} is to take advantage of the fact that the non-local operator $\mathcal{L}$ differs from the  fractional Laplacian by a zero order operator. Using this idea, we obtain a Maximum Principle for domains with small volume through the  ABP-estimate given  Proposition \ref{pro abp1} and we are able to use the moving planes method as in the case of the fractional Laplacian.
We prove first
\begin{proposition}\label{pro 7}
 Let ${\Sigma_\lambda}$ and $ {\Sigma_\lambda^-}$ be defined  as in the Section \S 3.
 Suppose that
$\varphi\in L^\infty({\Sigma_\lambda})$ and that $w_\lambda\in L^\infty(\R^N)\cap C(\R^N)$ is a solution of
\begin{equation}\label{eq a11}
\left\{ \arraycolsep=1pt
\begin{array}{lll}
-\mathcal{L} w_\lambda(x)\le \varphi(x)w_\lambda(x),\ \ \ \
&
x\in {\Sigma_\lambda},\\[2mm]
 w_\lambda(x)\geq0,& x\in \R^N\setminus {\Sigma_\lambda},
\end{array}
\right.
\end{equation}
where $\mathcal{L} $ was defined in (\ref{eq 2}).
Then, if $|{\Sigma_\lambda^-}|$ is small enough, $w_\lambda$ is non-negative in ${\Sigma_\lambda}$, that is,
$${w_\lambda}\ge 0\ \ \rm{in}\ \ {\Sigma_\lambda}.$$
 \end{proposition}
 {\bf Proof.} We define
$w^+_\lambda(x)$ as in \equ{eq c},
then we have
 \begin{eqnarray*}
 \mathcal{L}w^+_\lambda(x)
 &=&\int_{ B_1(x)}\frac{w_\lambda^+(x)-w_\lambda^+(z)}{|x-z|^{N+2\alpha_1}}dz+
\mu\int_{\R^N\setminus
B_1(x)}\frac{w_\lambda^+(x)-w_\lambda^+(z)}{|x-z|^{N+2\alpha_2}}dz
\\ &=&(-\Delta)^{\alpha_1}w^+_\lambda(x)\\
&&+\int_{\R^N\setminus
B_1(x)}(w_\lambda^+(x)-w_\lambda^+(z))(\frac\mu{|x-z|^{N+2\alpha_2}}-\frac1{|x-z|^{N+2\alpha_1}})dz
\\&\le&(-\Delta)^{\alpha_1}w^+_\lambda(x)+2C_0\|w^+_\lambda\|_{L^\infty({\Sigma_\lambda^-})}\ ,\ \ \ x\in{\Sigma_\lambda^-},
 \end{eqnarray*}
 where $C_0=\int_{\R^N\setminus B_1}|\frac\mu{|y|^{N+2\alpha_2}}-\frac1{|y|^{N+2\alpha_1}}|dy$. Thus we have
\begin{equation}\label{eq 122}
\Delta^{\alpha_1}w^+_\lambda(x)\leq -
\mathcal{L}w^+_\lambda(x)+2C_0\|w^+_\lambda\|_{L^\infty({\Sigma_\lambda^-})}\
,\ \ \ x\in{\Sigma_\lambda^-}.
\end{equation}
Since  $K_\mu$ is  radially
symmetric and decreasing in $|x|$, we may  repeat the arguments used to prove \equ{claim1} to get
\begin{equation}\label{eq 31}
\mathcal{L} w_\lambda^-(x)\le0,\ \ \ \ \forall\ x\in
\Sigma_\lambda^-,
\end{equation}
where $0<\lambda<1$ and  $w_\lambda^-$ was defined in (\ref{eq 4.01}).
Using (\ref{eq 122}), the linearity of $\mathcal{L}$,   (\ref{eq 31})  and equation (\ref{eq a11}), for all $x\in \Sigma_\lambda^-$,  we have
\begin{eqnarray}
\nonumber \Delta^{\alpha_1}w^+_\lambda(x) &\leq&
- \mathcal{L}w_\lambda(x)+\mathcal{L}w^-_\lambda(x)+2C_0\|w^+_\lambda\|_{L^\infty({\Sigma_\lambda^-})}
\\ \nonumber &\leq& - \mathcal{L}w_\lambda(x)+2C_0\|w^+_\lambda\|_{L^\infty({\Sigma_\lambda^-})}
\\ &\leq&\varphi(x)w_\lambda(x)+2C_0\|w^+_\lambda\|_{L^\infty({\Sigma_\lambda^-})}
 \leq C_1\|w^+_\lambda\|_{L^\infty({\Sigma_\lambda^-})},
 \end{eqnarray}
where $C_1=\|\varphi\|_{L^\infty({\Sigma_\lambda})}+2C_0$ and we notice that $w_\lambda=w_\lambda^+$ in $\Sigma_\lambda^-$.
Hence, we have
\begin{equation}\label{eq 124}
\left\{ \arraycolsep=1pt
\begin{array}{lll}
\Delta^{\alpha_1}w^+_\lambda(x)\le
C_1\|w^+_\lambda\|_{L^\infty({\Sigma_\lambda^-})},\ \ \ \ &
x\in {\Sigma_\lambda^-},\\[2mm]
 w^+_\lambda(x)=0,& x\in \R^N\setminus {\Sigma_\lambda^-}.
\end{array}
\right.
\end{equation}
Then, using Proposition \ref{pro abp1} with
$h(x)=C_1\|w^+_\lambda\|_{L^\infty({\Sigma_\lambda^-})}$, we obtain
a constant $C>0$  such that
 \begin{eqnarray*}
\|w^+_\lambda\|_{L^\infty(\Sigma_\lambda^-)} =
-\inf_{\Sigma_\lambda^-}w^+_\lambda \leq C d^{\alpha_1}
\|w^+_\lambda\|_{L^\infty(\Sigma_\lambda^-)}
|\Sigma_\lambda^-|^{\frac{\alpha_1} N},
\end{eqnarray*}
where $d=diam (\Sigma_\lambda^-)$. If $|\Sigma_\lambda^-|$ is small
enough we conclude that
$\|w_\lambda\|_{L^\infty(\Sigma_\lambda^-)}=\|w^+_\lambda\|_{L^\infty(\Sigma_\lambda^-)}=0,$
from where we
complete the proof.  \hfill$\Box$\\

\medskip
Now we provide a proof for  Theorem \ref{teo 5}.

\noindent{\bf Proof of Theorem \ref{teo 5}.} The proof of this theorem goes like the one for  Theorem \ref{teo 1} where
we use
Proposition \ref{pro 7}  instead of Proposition \ref{pro abp1}  and  $\mathcal{L}$ instead of  $(-\Delta)^{\alpha}$. The only
place where there is a difference is in the following property: for $0<\lambda<1$, if
 $w_\lambda\ge0$ and $w_\lambda\not\equiv0$
 in $\Sigma_\lambda$, then  $w_\lambda>0$
 in $\Sigma_\lambda$.

For $\mu\in(0,1]$, since $K_\mu$ is radially
symmetric and strictly decreasing, the proof of the property is similar to that given in Theorem \ref{teo 1}.
So we only need to prove it in case
$\mu=0$ so the kernel $K_0$ vanishes outside the unit ball $B_1$.
Let us assume that
 $w_\lambda\ge0$ and $w_\lambda\not\equiv0$
 in $\Sigma_\lambda$ and, by contradiction, let us  assume
 $\Sigma_0=\{x\in\Sigma_\lambda \ |\ w_\lambda(x)
=0\}\not=\O$.
By our assumptions on $w_\lambda$  we have that
$\Sigma_\lambda\setminus\Sigma_0=\{x\in\Sigma_\lambda \ |\
w_\lambda(x)>0\}$ is open and nonempty. Let us consider  $x_0\in \Sigma_0$
such that
\begin{equation}\label{eq 007}
dist(x_0,\Sigma_\lambda\setminus\Sigma_0)\le1/2,
\end{equation}
and observe that $(\Sigma_\lambda\setminus \Sigma_0)\cap B_1(x_0)$ is nonempty.
Using \equ{eq 51}  we have
\begin{eqnarray}
\mathcal{L} w_\lambda(x_0)&=& \mathcal{L}
u_\lambda(x_0)-\mathcal{L} u(x_0)\nonumber
\\&=&f(u_\lambda(x_0))-f( u(x_0))+g((x_0)_\lambda)-g(x_0)\nonumber
\\&=&g((x_0)_\lambda)-g(x_0)\ge 0,\label{eq 07}
\end{eqnarray}
where the last inequality holds by monotonicity assumption on $g$ and since
   $|x_0|>|(x_0)_\lambda|$.
On the other hand, denoting by $A_\lambda=\{(x_1,x')\in\R^N\ |\
x_1>\lambda\}$, since $w_\lambda(x_0)=0$ and $w_\lambda(z_\lambda)=-w_\lambda(z)$ for any
$z\in\R^N$,  we have
\begin{eqnarray*}
\mathcal{L} w_\lambda(x_0)
&=&-\int_{A_\lambda}w_\lambda(z){K_0}(x_0-z)dz
-\int_{\R^N\setminus A_\lambda}w_\lambda(z){K_0}(x_0-z)dz
\\&=&-\int_{A_\lambda}w_\lambda(z){K_0}(x_0-z)dz-\int_{A_\lambda}w_\lambda(z_\lambda){K_0}(x_0-z_\lambda)dz
\\&=&-\int_{A_\lambda}w_\lambda(z)({K_0}(x_0-z)-{K_0}(x_0-z_\lambda))dz.
\end{eqnarray*}
Since $|x_0-z_\lambda|>|x_0-z|$ for $z\in A_\lambda$, by definition of $K_0$, $\Sigma_\lambda$ and $\Sigma_0$, we have that
$${K_0}(x_0-z)>{K_0}(x_0-z_\lambda)\quad\mbox{and}\quad    w_\lambda(z)>0   \quad  {\rm{for}}\quad z\in(\Sigma_\lambda\setminus \Sigma_0)\cap B_1(x_0),$$
and we also have that
$w_\lambda(z)\geq0$ and ${K_0}(x_0-z)\geq{K_0}(x_0-z_\lambda)$ for all $z\in A_\lambda,$
so that
$$\mathcal{L} w_\lambda(x_0)<0,$$
contradicting  (\ref{eq 07}).
Hence $\Sigma_0$ is empty and then
$w_\lambda>0$ in $\Sigma_\lambda$, completing
 the proof of the theorem.
 \hfill$\Box$\\

\begin{remark}\label{remark add2}
The theorem we just proved  can be extended to more general non-homogeneous kernels in the following class
\begin{equation}\label{eq add120}
K(x)=\left\{ \arraycolsep=1pt
\begin{array}{lll}
 |x|^{-N-2\alpha},\ \ \ \ & x\in B_{r},\\[2mm]
 \theta(x),\ \ \ \ & x\in B^c_{r},
\end{array}
\right.
\end{equation}
here $\alpha\in(0,1)$, $r>0$ and the function $\theta:B^c_{r}\to\R$ satisfies that
\begin{itemize}
\item[$(C)\ $]
$\theta\in L^1(B^c_{r})$ is nonnegative, radially symmetric and such that the kernel $K$ is decreasing.
\end{itemize}
\end{remark}

\noindent {\bf Acknowledgements:}  P.F. was  partially supported by Fondecyt Grant \# 1110291,
 BASAL-CMM projects and CAPDE, Anillo ACT-125.
Y.W. was partially supported by Becas CMM.

\end{document}